\documentclass{amsart}

\usepackage{amscd}
\usepackage{amssymb}
\usepackage[all]{xy}
\usepackage{mathrsfs}
\numberwithin{equation}{subsection}

\newtheorem{theorem}{Theorem}[subsection]
\newtheorem{corollary}[theorem]{Corollary}
\newtheorem{lemma}[theorem]{Lemma}
\newtheorem{proposition}[theorem]{Proposition}

\theoremstyle{definition}
\newtheorem{defn}[theorem]{Definition}
\newtheorem{example}[theorem]{Example}
\newtheorem{remark}[theorem]{Remark}

\newcommand\gfrac[2]{\genfrac{}{}{0pt}{}{#1}{#2}}
\newcommand\isomarrow{\stackrel{\cong}{\longrightarrow}}

\title[Alcove walks and nearby cycles on affine flag manifolds]{Alcove walks\\and nearby cycles on affine flag manifolds}
\author{Ulrich G\"ortz}
\address{Mathematisches Institut, Beringstr. 1, 53115 Bonn, Germany}
\email{ugoertz@math.uni-bonn.de}
\date{}

\begin{document}

\maketitle

\begin{abstract}
Using Ram's theory of alcove walks, we give a proof of the Bernstein
presentation of the affine Hecke algebra. The method works also in the case
of unequal parameters. We also discuss how these results help in studying
sheaves of nearby cycles on affine flag manifolds.
\end{abstract}
\section{Introduction}

\subsection{}
In a recent paper \cite{R}, Ram has introduced the notion of \emph{alcove
walk} and used it in order to describe the affine Hecke algebra associated
to a root datum. 

In these notes we will show that, with a little extra work, this concept
yields a proof of the Bernstein presentation of the affine Hecke algebra
$\mathcal H$. The method applies to the case of unequal parameters as well
(see section \ref{unequalparams}), and we obtain a proof which might be
considered less technical than the one given by Lusztig \cite{L}.

The main new ingredient beyond Ram's paper is the following
theorem which we will state here in the introduction without using the
language of alcove walks. Let $W_a$ be the affine Weyl group associated
with some root system, and let $s_0, \dots, s_r$ denote the simple
reflections, which generate $W_a$. Denote by $\mathbf a$ the base alcove
in the ``standard apartment''.

\begin{theorem} Let $w\in W_a$. For an expression
\begin{equation} \label{expr}  w = s_{i_1} \cdots s_{i_n}
\end{equation}
of $w$ as a product in the generators (which does not have to be reduced),
consider the element 
\[ \Psi(w) := T_{s_{i_1}}^{\varepsilon_1} \cdots T_{s_{i_n}}^{\varepsilon_n} \]
in the affine Hecke algebra, where the $\varepsilon_\nu \in\{ \pm 1 \}$ are
determined as follows. Let $\mathbf b$ be an alcove far out in the
anti-dominant chamber (``far out'' depends on $w$, and the result will then
be independent of $\mathbf b$, see section \ref{orientations} for a
precise definition). For each $\nu$, consider the alcove $\mathbf c_\nu
:= s_{i_1}\cdots s_{i_{\nu-1}} \mathbf a$, and denote by $H_\nu$ the
affine root hyperplane containing its face of type $i_\nu$.
We set
\[ \varepsilon_\nu = \left\{ \begin{array}{ll} 
1 & \text{if } \mathbf c_\nu \text{ is on the
same side of } H_\nu \text{ as } \mathbf b \\
-1 & \text{otherwise} 
\end{array} \right.\]
Then the element $\Psi(w)$ is independent of the choice of expression
(\ref{expr}).
\end{theorem}
This is Theorem \ref{indep} in the text. It is clearly implicitly
contained in Ram's paper,
and I actually learned its statement from Ram before his paper was
finished. However, I do not quite see how to extract the theorem from
what is spelled out in \cite{R} (cf. remark \ref{gap}). 
In any case, it may be interesting to have
the proof of the theorem below, since it can then be used to infer, from
Ram's arguments, a proof of the Bernstein presentation of the affine
Hecke algebra, starting from the Iwahori-Matsumoto presentation. 
Ram takes the Bernstein presentation as the definition of
the affine Hecke algebra.

Another application is the existence of so-called minimal expressions (cf.
the paper \cite{HP} by Haines and Pettet) for the elements $\Theta_\lambda
\in \mathcal H$.

\subsection{Acknowledgments} I am grateful to Arun Ram for explaining to me his
theory of alcove walks, and in particular for showing me the statement of
Theorem \ref{indep}. I also thank Thomas Haines who pointed out a few
inaccuracies in a preliminary version of this text and made several very
helpful remarks.

\section{The alcove walk algebra}
\subsection{Notation}

In this section we will collect the relevant notation; for more details on
these notions, see Humphreys' book \cite{Hu} and the papers by Lusztig
\cite{L}, Haines and Pettet \cite{HP} and Haines, Kottwitz, and Prasad
\cite{HKP}, for instance. Note that Lusztig uses a setup which is
dual to ours: he works with roots where we use coroots, and conversely.

Let $(X^*, X_*, R, R^\vee, \Delta)$ be a reduced and irreducible based root
datum with $\Delta$ being the set of simple roots, and denote by $W$ its
Weyl group, generated by the simple reflections $\{ s_\alpha;\ \alpha \in
\Delta \}$. Denote by $\widetilde{W} := X_* \rtimes W$ the extended affine
Weyl group. For $\lambda\in X_*$, we denote by $\epsilon^\lambda$ the
corresponding element in $\widetilde{W}$.  Let $S_a = \{ s_\alpha;\ \alpha
\in \Delta\} \cup \{ s_0 \}$, where $s_0 =
\epsilon^{\alpha_0^\vee}s_{\alpha_0^\vee}$ and where $\alpha_0$ is the
unique highest root. The subgroup $W_a \subseteq \widetilde{W}$ generated
by $S_a$ is the affine Weyl group of the root system associated with our
root datum, and $(W_a, S_a)$ is a Coxeter system.

We define a length function $\ell \colon \widetilde{W} \longrightarrow
\mathbb Z$ as follows:
\[ \ell(w \epsilon^\lambda) = \sum_{\gfrac{\alpha>0}{w(\alpha) < 0}} |
\langle \alpha, \lambda \rangle +1 | + \sum_{\gfrac{\alpha>0}{w(\alpha)>0}}
| \langle \alpha,\lambda \rangle |. \]
This function extends the length function on $W_a$.
We have a short exact sequence
\[ 1 \longrightarrow  W_a \longrightarrow \widetilde{W} \longrightarrow X_*/Q^\vee
\longrightarrow 0, \]
where $Q^\vee$ is the coroot lattice, i.~e.~the subgroup of $X_*$ generated
by $R^\vee$. The restriction of the projection $\widetilde{W}
\longrightarrow X_*/Q^\vee$ to the subgroup $\Omega \subseteq
\widetilde{W}$ of elements of length $0$ is an isomorphism $\Omega
\isomarrow X_*/Q^\vee$.

We extend the Bruhat order on $W_a$ by declaring 
\[ w\tau \le w'\tau' \Longleftrightarrow w \le w', \tau=\tau',\qquad
w,w'\in W_a,\ \tau, \tau' \in \Omega. \]

For an affine root $\beta = \alpha - n$, $\alpha\in R$, $n \in \mathbb Z$,
we have the hyperplane $H_\beta = H_{\alpha,n} = \{ x \in X_{*,\mathbb R};\
\langle \alpha, x \rangle = n \}$.  An alcove is a connected component of
the complement of the union of all affine root hyperplanes inside
$X_{*,\mathbb R}$. The choice of $\Delta$ determines a base alcove $\mathbf
a$ which is the unique alcove contained in the dominant finite Weyl chamber
whose closure contains the origin.  The group $\widetilde{W}$ acts on
$X_{*,\mathbb R}$, and since the union of all affine root hyperplanes is
stable under this action, we have an action of $\widetilde{W}$ on the set
of alcoves. The affine Weyl group $W_a$ acts simply transitively on the set
of alcoves, so we can identify $W_a$ with the set of alcoves in the
standard apartment $X_{*,\mathbb R}$ by mapping $w \in W_a$ to the alcove
$w\mathbf a$. On the other hand, the group $\Omega$ of elements of length
$0$ in $\widetilde{W}$ is precisely the stabilizer of the base alcove
inside $\widetilde{W}$. If $\lambda$ denotes the image of the origin under
$\tau \in \Omega$, then $\tau w_0 w = \varepsilon^\lambda$, where $w$ is
the longest element in the stabilizer $W_\lambda$ of $\lambda$ in $W$, and
$w_0$ is the longest element in $W$.

Let $r = \# \Delta$ be the semi-simple rank of the root datum. We
order the simple reflections in some way and denote them by $s_1, \dots,
s_r$ in the sequel.
The group $\Omega$ acts on the set of simple affine reflections (resp.~on
the set of simple affine roots) and for $\tau \in \Omega$ we define
$\tau(i)$ by $\tau s_i \tau^{-1} = s_{\tau(i)}$.

Let us also briefly recall the definition of the affine Hecke algebra.
We fix a ground ring $k$, an invertible element $v \in k$ and set $q =
v^2$. For example, $k$ might be a field (as in \cite{R}), or we could let
$k = \mathbb Z[v, v^{-1}]$ be the ring of Laurent polynomials over the integers.
We will often write $q^{\frac{1}{2}}$ instead of $v$.
The braid group of $\widetilde{W}$ is the group with generators
\[ T_w, \quad w \in \widetilde{W}, \]
and relations
\[ T_w T_{w'} = T_{ww'} \text{ for } w, w' \in \widetilde{W} \text{ with }
\ell(ww') = \ell(w) + \ell(w'). \]
The affine Hecke algebra $\mathcal H$ is the quotient of the group algebra
of the braid group (over our
fixed ground ring $k$), by the two-sided ideal generated by
\[ (T_s + 1)(T_s-q),\quad s \in S_a. \]
We denote the image of $T_w$ in $\mathcal H$ again by $T_w$.
We sometimes abbreviate $T_{s_i}$ to $T_i$. Further, it is often useful to
use the element $\tilde{T}_w := q^{-\ell(w)/2}T_w$ instead of $T_w$. We
use $\tilde{T}_i$ as an abbreviation for $\tilde{T}_{s_i}$.

For $\lambda\in X_*$, we define $\Theta_\lambda\in\mathcal H$ as follows:
write $\lambda = \lambda_1 - \lambda_2$ with $\lambda_1$, $\lambda_2$
dominant, and let $\Theta_\lambda = \tilde{T}_{\epsilon^{\lambda_1}}
\tilde{T}_{\epsilon^{\lambda_2}}^{-1}$ (which is well-defined as an element
of $\mathcal H$).

It is not hard to see that the elements $T_w$, $w\in \widetilde{W}$ form a
$k$-basis of $\mathcal H$. We note the following lemma which exhibits
variants of this basis and which is easily proved by using that $T_i^{-1}
= q^{-1}T_i + (q^{-1}-1)$.

\begin{lemma} \label{bases}
Fix a reduced expression $w = s_{i_1} \cdots s_{i_\ell}\tau$ for each $w\in
\widetilde{W}$, and fix signs $\varepsilon_i(w)\in \{\pm 1\}$. Then the set
of all elements $T_{i_1}^{\varepsilon_1(w)}\cdots
T_{i_\ell}^{\varepsilon_\ell(w)}T_\tau \in \mathcal H$ is a $k$-basis of
$\mathcal H$.
\end{lemma}

Similarly, the elements $\Theta_\lambda T_w$,
$\lambda \in X_*$, $w\in W$ (and likewise the elements $T_w
\Theta_\lambda$) form a $k$-basis of $\mathcal H$. (See \cite{L},
Prop.~3.7, or \cite{HKP} Lemma 1.7.1.)


\subsection{Definition of the alcove walk algebra}

The alcove walk algebra $A$ is the (non-commutative) $k$-algebra with
generators 
\[ c_i^+, c_i^-, f_i^+, f_i^-, \quad i = 0,\dots, r, \quad t_\tau,\ \tau \in \Omega \]
and relations 
\begin{eqnarray*}
&& c_i^- = c_i^+ + f_i^-, \quad f_i^- = - f_i^+, \quad i = 0,\dots, r, \\
&& t_\tau ?_i = ?_{\tau(i)} t_\tau, \quad ?\in\{c^+, c^-, f^+, f^-\},\ \tau \in
\Omega, \\
&& t_\tau t_{\tau'} = t_{\tau\tau'}, \quad \tau, \tau' \in \Omega
\end{eqnarray*}
The elements $c_i^+$, $c_i^-$, $f_i^+$, $f_i^-$ are called the positive
crossing, the negative crossing, the positive folding and the negative
folding of type $i$, respectively.  We have a natural map $\Phi\colon A
\longrightarrow \mathcal H$ from the alcove walk algebra to the affine
Hecke algebra by mapping
\[ c_i^+ \mapsto \tilde{T}_i,\quad  c_i^- \mapsto \tilde{T}_i^{-1},\quad
f_i^+ \mapsto q^{\frac{1}{2}}-q^{-\frac{1}{2}},\quad f_i^- \mapsto
q^{-\frac{1}{2}}-q^{\frac{1}{2}},\quad t_\tau\mapsto T_\tau.\] 
We will determine the kernel of this map (see Proposition \ref{kernel}),
and hence get a new description of the affine Hecke algebra as a certain
quotient of the alcove walk algebra.

\subsection{Orientation}
\label{orientations}

In order to explain the terminology \emph{alcove walk} algebra, we
fix an orientation in the following sense.

\begin{defn}
A \emph{root hyperplane orientation} is given by distinguishing, for each
affine root hyperplane $H$, a positive half-space among
the two half-spaces which form the complement of $H$ in
$X_{*,\mathbb R}$, such that either 
\begin{enumerate}
\item for any finite set of affine root hyperplanes, the intersection of
the corresponding negative half-spaces is non-empty (and hence contains an
alcove), or
\item for any finite set of affine root hyperplanes, the intersection of
the corresponding positive half-spaces is non-empty
\end{enumerate}
(Here the half-space in $X_{*,\mathbb R}\setminus H$ which is not
positive, is called negative.)
\end{defn}
Given an orientation of type (1), and an alcove $\mathbf b$ in the
intersection of the negative half-spaces associated to a 
fixed finite set of affine root hyperplanes, we can say
that for these hyperplanes the orientation is given by prescribing that
crossing the hyperplane ``in the positive direction'', i.~e.~from the
negative to the positive half-space, is the same as ``walking away'' from
$\mathbf b$.

There are two obvious examples for orientations:
\begin{example}
If $\mathbf b$ is a fixed alcove, we can just use it to define an
orientation by saying that for each hyperplane, the negative half-space is
the one containing $\mathbf b$. We can express this by saying that the
most negative point of the orientation is lying inside $\mathbf b$.
Analogously, we get another orientation by saying that for each hyperplane,
the positive half-space is the one containing $\mathbf b$.
\end{example}

\begin{example}
The orientation which will be most important for us is given by calling,
for each \emph{positive} root $\alpha$, and each $j\in \mathbb Z$, the half-space
\[ \{ x \in X_{*,\mathbb R};\ \langle x, \alpha \rangle > j \} \] 
the positive half-space. In other words (cf \cite{BT}), the negative
half-space is the unique half-space which contains a quartier of the form
$y + C^-$, where $C^-$ denotes the anti-dominant Weyl chamber. We can
describe this orientation by saying that the most negative point lies
infinitely deep in the anti-dominant chamber.  We will call this
orientation the \emph{standard orientation}; it is the one used in
\cite{R}.  It also occurs in the paper \cite{GL} by Gaussent and
Littelmann.

Alternatively, we could replace the anti-dominant chamber by any other Weyl
chamber. Orientations of this type implicitly play a role in \cite{GHKR}.
\end{example}

\subsection{Alcove walks}
We will now give a formal definition of the notion of \emph{alcove walk},
and simultaneously define the end-point ${\rm end}(\gamma)\in
\widetilde{W}$ for an alcove walk $\gamma$.
For a more ``pictorial'', and probably more accessible
definition, see \cite{R}.

\begin{defn} 
Fix an orientation as defined in the previous section.  Alcove walks are
pairs $(w, \tau)$ where $w$ is a word in the $c_i^+$, $c_i^-$, $f_i^+$,
$f_i^-$, $i\in\{0,\dots, r\}$, and where $\tau \in \Omega$, subject to
certain conditions. The length of the alcove walk is by definition the
length of the word $w$. We define the conditions an alcove walk has to
satisfy by induction on its length.
\begin{enumerate}
\item If $\tau\in\Omega$, then $(\emptyset, \tau)$ is an alcove walk of
length $0$ (where $\emptyset$ denotes the empty word). Its end point is
${\rm end}((\emptyset, \tau)) = \tau \in \widetilde{W}$.
\item If $(w, \tau)$ is an alcove walk, with $w = w_1 \cdots w_{n}$ a
word of length $n$, and $w_{n+1} \in \{c_i^+, c_i^-, f_i^+, f_i^-\}$, such
that either
\begin{itemize}
\item ${\rm end}((w, 0))\mathbf a$ is on the negative side of the wall of type $i$
adjacent to the alcove ${\rm end}((w,0))\mathbf a$, and $w_{n+1} \in \{
c_i^+, f_i^-\}$, or
\item ${\rm end}((w, 0))\mathbf a$ is on the positive side of the wall of type $i$
adjacent to the alcove ${\rm end}((w,0))\mathbf a$, and $w_{n+1} \in \{
c_i^-, f_i^+ \}$,
\end{itemize}
then $(ww_{n+1}, \tau)$ is an alcove walk, and its end point is
\[ {\rm end}((ww_{n+1}, \tau)) = \left\{ \begin{array}{ll} {\rm
end}((w,0))s_i\tau & \text{if } w_{n+1} \in \{ c_i^+, c_i^- \} \\ {\rm
end}((w,\tau)) & \text{if } w_{n+1} \text{ has the form } f_i^\pm.
\end{array} \right.\]
\item All alcove walks arise in this way.
\end{enumerate}
\end{defn}

Since we put the $\Omega$-part into the second component, the property of
being an alcove walk is actually independent of the $\Omega$-component, and
furthermore for the end points we have ${\rm end}((w,\tau)) = {\rm
end}((w,0))\tau$. Because in the definition given here, we build the walks
``from left to right'', we always have to insert the next step between the
given walk and the $\tau$, so that the relevant information about the
orientation of the adjacent walls is given by ${\rm end}((w,0))$ rather
than by ${\rm end}((w,\tau))$.

The sequence of end points 
\[{\rm end}((w_1, 0)), {\rm end}((w_1w_2, 0)),
\dots, {\rm end}((w_1\dots w_n, 0)) \in W_a \]
should be seen as a
sequence of alcoves in the standard apartment $X_{*,\mathbb R}$---hence the
name \emph{alcove walk}.

We call an alcove walk $(w, \tau)$ \emph{non-folded}, if no symbols of the
form $f_i^+, f_i^-$ occur in the word $w$. We consider alcove walks as
elements of the alcove walk algebra in the obvious way. 
As the following lemma shows, we
can see the choice of orientation as the choice of a basis of the alcove
walk algebra.

\begin{lemma}
The set of alcove walks is a basis of the alcove walk algebra as a
$k$-module.
\end{lemma}

\begin{proof}
This is \cite{R}, Lemma 3.1; since Ram does not give a proof, for the
convenience of the reader we produce a proof here. To simplify the
notation, let us suppose that $\Omega = \{ 0 \}$. Because of the relations
$c_i^- = c_i^+ + f_i^-$, $f_i^+ = -f_i^-$, it is clear that $A$ is
isomorphic to the free algebra with generators $c_i^+$, $f_i^-$. Hence the
set $\mathscr B$ of words in $\{c_i^+, f_i^-;\ i=0,\dots, r \}$ is a 
$k$-basis of $A$.
Now fix an orientation and denote by $\mathscr W$ the set of alcove walks.
We have an obvious bijection $\mathscr W \longrightarrow \mathscr B$ by
mapping each walk to the element of $\mathscr B$ by changing all exponents
of $c$'s to $+$, and all exponents of $f$'s to $-$. We order $\mathscr
W$ in some way, such that whenever a walk $w_1$ has more $c$'s in it than a
walk $w_2$, then $w_1 > w_2$. The bijection $\mathscr W \cong \mathscr B$
induces an order on $\mathscr B$ with the same property. Now if we express
each walk in $\mathscr W$ in terms of the basis $\mathscr B$ and take the
coefficients as the column vector of an (infinite) square matrix (using the order we
defined), then this matrix will be upper triangular with $1$'s on the
diagonal, and it follows that $\mathscr W$ is a basis as well.
\end{proof}

\section{The independence result}

\subsection{}
We again fix an orientation as defined in section \ref{orientations}.
Given a word $s_{i_1} \cdots s_{i_k}\tau$ (not necessarily reduced) in the
extended affine Weyl group, we can associate to it, or to the corresponding
gallery, a unique non-folded alcove walk $c_{i_1}^{\varepsilon_1} \cdots
c_{i_k}^{\varepsilon_k} t_\tau$, $\varepsilon_i \in \{ +, -\}$. On the other
hand, we can associate to the given word the element
\[T_{i_1}^{\varepsilon_1} \cdots T_{i_k}^{\varepsilon_k} T_\tau \]
in the affine Hecke algebra. We also denote this element by
\[T_{i_1}^{\varepsilon} \cdots T_{i_k}^\varepsilon \tau,\]
i.~e.~we let $\varepsilon$ denote the appropriate sign, depending on $\nu
\in \{ 1, \dots, k \}$.

As a variant, we can consider alcove walks which do not start at the base
alcove, but at another alcove, say at $w\mathbf a$. We will denote by
\[T_{i_1}^{\varepsilon(w)} \cdots T_{i_k}^{\varepsilon(w)}\]
the corresponding element in $\mathcal H$, where again $\varepsilon(w)$ is
understood to vary with $\nu \in \{ 1, \dots, k \}$.

\begin{theorem} \label{indep}
Let $w \in \widetilde{W}$ be an element in the extended affine Weyl group,
and let 
\[ w = s_{i_1} \cdots s_{i_k}\tau = s_{j_1} \cdots s_{j_\ell}\tau\]
be expressions for $w$ (which need not be reduced). Then in the
affine Hecke algebra, we have the equality
\[ T_{i_1}^{\varepsilon} \cdots T_{i_k}^\varepsilon T_\tau = T_{j_1}^\varepsilon
\cdots T_{j_\ell}^\varepsilon T_\tau.\]
\end{theorem}

\begin{proof}
We clearly may assume that $\tau={\rm id}$, i.~e.~that $w$ actually is an
element of the affine Weyl group $W_a$. 
The affine Weyl group is a Coxeter group, and hence has the so-called
\emph{word property} (see \cite{T}, \cite{BB}, Theorem 3.3.1), i.~e.~we can
get the expression $s_{j_1} \cdots s_{j_\ell}$ from $s_{i_1} \cdots
s_{i_k}$ by applying transformations of the following kinds (in a suitable
order):
\begin{enumerate}
\item nil-move: delete a subexpression of the form $s_i s_i$ from the word
\item inverse nil-move: insert a subexpression of the form $s_is_i$
somewhere
\item braid move: replace a subexpression  $s_i s_{i'} s_i \cdots $ by
$s_{i'} s_i s_{i'} \cdots$, where both these words consist of $m_{i,i'}$
letters, $m_{i,i'}$ being the entry in the Coxeter matrix corresponding to
$i$, $i'$.
\end{enumerate}
Hence it is enough to prove that these types of transformations do not
change the product $T_{i_1}^{\varepsilon} \cdots T_{i_k}^\varepsilon$ (of
course the signs $\cdot^\varepsilon$ have to be taken into account). For
nil-moves and inverse nil-moves this is obvious, since the two adjacent
$T$'s will have exponents $1$ and $-1$ (or $-1$ and $1$), thus will cancel.

What remains to show is that for all $i,j$, and for all $w\in W_a$, we have
\begin{equation} \label{eq1}
T_i^{\varepsilon(w)} T_j^{\varepsilon(w)} T_i^{\varepsilon(w)} \cdots =
T_j^{\varepsilon(w)} T_i^{\varepsilon(w)} T_j^{\varepsilon(w)} \cdots,
\qquad (*)
\end{equation}
where both products have $m_{i,j}$ factors.
(Of course, in the case that all the $\varepsilon(w)$'s are equal, the
equality follows immediately from the braid relations in the affine Weyl
group.)

Since only finitely many alcoves are involved in the alcove walk, by the
definition of orientation, there is an alcove $\mathbf b= v\mathbf a$, $v
\in W_a$, such that the positive/negative direction is determined by
whether we are approaching $\mathbf b$, or not. We may assume without loss
of generality that $\mathbf b$ is the most negative point (rather than the
most positive point) for the finitely many hyperplanes involved, because
otherwise we could replace the orientation by its ``inverse'': the equality
we have to check is the same for both of these orientations.
Denote by $\mathbf o$ the orientation given by making the base alcove the
most negative point, and by $\varepsilon_{\mathbf o}$ the signs defined
with respect to $\mathbf o$.  We then have
\[\varepsilon(w) = \varepsilon_{\mathbf o}(v^{-1}w) \]
in the obvious sense.  Hence it is enough to check the assertion of the
theorem for the orientation $\mathbf o$, i.~e. we may assume that the signs
are determined by whether we come closer to a fixed point in the interior
of the base alcove, or not; in other words whether the length of the
element in the affine Weyl group corresponding to the alcove decreases, or
increases.

The coset $w W_{i,j}$ of the parabolic subgroup $W_{i,j} \subseteq W_a$
generated by $s_i$, $s_j$ has a unique element of minimal length, and a
unique element of maximal length, and there are two ways to go from the
minimal length element to the maximal length element.  We multiply
(\ref{eq1}) on the right 
by the inverse of the right hand side, and get an equation of the form 
\[ T_i^{\varepsilon_1} T_j^{\varepsilon_2} \cdots
T_j^{\varepsilon_{2m_{i,j}}} = 1,\] 
and the alcove walk corresponding to the left hand side of this equation
starts at $w\mathbf a$ and then comes back to $w\mathbf a$, seeing each
alcove in the coset $w W_{i,j}$ exactly once. Since the situation is
symmetric, we may and will assume that we start in the positive direction.
Then we continue in the positive direction until we get to the element of
maximal length, from there we go in the negative direction until we get to
the element of minimal length, and finally go in the positive direction
back to $w\mathbf a$, where we started.

The string of $T$'s corresponding to the part of the walk going from the
maximal to the minimal element is a string of length $m_{i,j}$, and all
$T$'s have an exponent $-1$. We may hence replace this string with the
string which has $i$ and $j$ exchanged because of the braid relations in
the affine Weyl group, and we then get that the whole product cancels. 
\end{proof}

Note that the statement of the proposition remains true, if we replace
$T_i$ by $\tilde{T}_i$ everywhere (the same proof applies).

\subsection{The kernel of $\Phi$}

As before, we fix an orientation. Recall that the notion of alcove walk
depends on the orientation, but that neither the definition of the alcove
walk algebra, nor the morphism $\Phi \colon A \longrightarrow \mathcal H$
do.  Therefore the following proposition is a little surprising.
\begin{proposition} \label{kernel}
The kernel of $\Phi\colon A\longrightarrow \mathcal H$ is the two-sided
ideal $\mathscr J \subseteq A$ generated by
\begin{eqnarray*}
&&  f_i^+ - (q^{\frac{1}{2}}-q^{-\frac{1}{2}}),\quad i=0,\dots, r \\
&& p - p' \text{ for } p, p' \text{ non-folded alcove walks with the same
end point}
\end{eqnarray*}
\end{proposition}
We remark that this definition implies in particular that $c_i^- -
(c_i^+)^{-1}\in \mathscr J$ for all $i$ (of course, these elements
obviously lie in $\ker \Phi$).
 
\begin{proof} Theorem \ref{indep} shows that $\Phi$ factors through a
morphism $A/\mathscr J \longrightarrow \mathcal H$ which of course is again
surjective. Fixing a non-folded walk $p_w$ from $\mathbf a$ to
$w\mathbf a$ for each $w\in \widetilde{W}$ gives us a set of elements in
$A/\mathscr J$ which generates $A/\mathscr J$ as a $k$-module. Lemma
\ref{bases} implies that this set is mapped to a basis of $\mathcal H$, and
hence the morphism $A/\mathscr J \longrightarrow \mathcal H$ is an
isomorphism.  \end{proof}

Each choice of orientation hence gives us a basis of $\mathcal H$
consisting of the images of non-folded walks to $w\mathbf a$, $w\in
\widetilde{W}$, so in a sense the choice of orientation corresponds to the
choice of a basis for $\mathcal H$; cf. Remark 3.6 in \cite{R}. Of course
we do not get every basis of $\mathcal H$ in this way.

It is easy to see (and not surprising) that neither Theorem
\ref{indep} nor the proposition above hold for ``orientations'' which do
not satisfy the condition imposed in section \ref{orientations}.

\begin{remark}
It seems that the notion of alcove walk is related to the Bruhat-Tits
building. Assume that $q$ is the number of
elements of the residue class field of a local field $K$, and fix a split
reductive algebraic group $G$ over $K$ which gives rise to the root system
under consideration.


Denote by $\rho$ the retraction of the Bruhat-Tits building to the standard
apartment from an alcove ``far out'' in the anti-dominant chamber. More
precisely, for each alcove in the building, its image under such
retractions will depend on the alcove, but will stabilize if the alcove is
sufficiently deep in the anti-dominant chamber, and this gives us the image
of the alcove under $\rho$. See \cite{BT} 2.9.1.

Then to each non-stuttering gallery in the building (starting at the base
alcove) we can associate via this retraction a unique alcove walk.
Let $\mathbf b$ be an alcove in the building which is part of such a
gallery. If $\rho(\mathbf b)$ is on the positive side of the wall of type $i$
adjacent to it, then there are $q-1$ alcoves $\mathbf b'$
adjacent to $\mathbf b$, but different from it, such that $\rho(\mathbf
b')=\rho(\mathbf b)$, and the alcove walk will have a positive folding
precisely if one of these alcoves $\mathbf b'$ is the successor of $\mathbf
b$ in the gallery. It may be possible to heuristically explain the $f_i^+$ in this
way (using a different normalization, such that $f_i^+$ will correspond to
$q-1$), but it is not clear to me how to establish an explicit and precise
relationship.

\end{remark}

\subsection{The Bernstein relations}

In this section we work with the standard orientation, i.~e.~we put the
most negative point infinitely deep in the anti-dominant chamber.
We define elements $t_w, \theta_\lambda \in \mathcal H$ ($w\in W$,
$\lambda\in X_*$) as follows.
For $w \in W$, let $p\in A$ be a non-folded alcove walk from $\mathbf a$ to
$w^{-1} \mathbf a$, and let $t_w := \Phi(p)^{-1} \in \mathcal H$.

For $\lambda \in X_*$, denote by $\theta_\lambda\in\mathcal H$ the image
under $\Phi$ of a non-folded alcove walk with end point $\epsilon^\lambda$.

The following proposition is Proposition 3.2 in \cite{R}. It shows that the
Bernstein relations are satisfied for the elements $\theta_\lambda$, $t_w$
in $\widetilde{H}$.

\begin{proposition} \label{Ram-Bernstein}
Let $\tau \in \Omega$, $\lambda, \mu \in X_*$, $w \in W$, and $1 \le i \le
r$, and denote by $\alpha_i$ the corresponding simple root. 
Let $\varphi \in R$ be the positive root such that $H_{\alpha_0} =
H_{\varphi,1} := \{ x \in X_{*,\mathbb R};\ \langle \varphi, x \rangle = 1
\}$ is the wall of $\mathbf a$ which is not a wall of the Weyl chamber
$\mathbf a$ lies in. Then we have
\begin{enumerate}
\item $\theta_\lambda \theta_\mu = \theta_{\lambda+\mu}$
\item $t_{s_i} t_w = \left\{ \begin{array}{ll}
t_{s_iw} & \text{if } \ell(s_iw) > \ell(w) \\
t_{s_iw} + (q^{\frac{1}{2}}-q^{-\frac{1}{2}})t_w & \text{otherwise}
\end{array} \right.$
\item $t_{s_i} \theta_\lambda = \theta_{s_i\lambda}t_{s_i} +
(q^{\frac{1}{2}}-q^{-\frac{1}{2}}) \frac{\theta_\lambda -
\theta_{s_i\lambda}}{1-\theta_{-\alpha_i^\vee}}$
\item $\Phi(c_0^+) t_{s_\varphi} = \theta_{\varphi^\vee}$, where $s_\varphi$ denotes
the reflection associated with $\varphi$.
\item Let $\tau\in \Omega$ be an element of $\widetilde{W}$ of length $0$.
Recall that $\tau$, as an automorphism of $X_{*,\mathbb R}$ maps the base
alcove to itself. Denote by $\lambda \in X_*$ the image of the origin under
$\tau$. Then
\[ \theta_\lambda = T_\tau t_{w_0w}, \]
where $w$ is the longest element in the stabilizer $W_\lambda$ of $\lambda$
in $W$, and $w_0$ is the longest element in $W$.
\end{enumerate}
\end{proposition}

\begin{proof} All these relations can be checked with relatively little
effort in terms of alcove walks. We give some of the details, since the
proof in \cite{R} is partly quite terse. If $p_\lambda$ is a
non-folded walk from $\mathbf a$ to $\epsilon^\lambda\mathbf a$, and
$p_\mu$ is a non-folded walk from $\mathbf a$ to $\epsilon^\mu\mathbf
a$, then clearly the composition $p_\lambda p_\mu$ is a non-folded walk
from $\mathbf a$ to $\epsilon^{\lambda+\mu}\mathbf a$, and this gives
(1). 

To get (3), we may assume without loss of generality that $\langle
\alpha_i, \lambda \rangle \ge 0$. Let $p_\lambda = (c_{i_1}^{\varepsilon_1}
\cdots c_{i_\ell}^{\varepsilon_\ell}, \tau)$ be a walk from $\mathbf a$ to
$\epsilon^\lambda \mathbf a$ of minimal length. Let $\mathbf a_1 = \mathbf
a$, $\mathbf a_2 = s_{i_1}\mathbf a$, $\dots$, $\mathbf a_\ell = s_{i_1}
\cdots s_{i_\ell} \mathbf a$ be the sequence of alcoves ``visited'' by this
path. Consider the element $c_i^- p_\lambda \in \mathcal H$; the
corresponding sequence of alcoves is the base alcove plus the mirror image
of the sequence $\mathbf a_1, \dots, \mathbf a_\ell$ with respect to the
reflection $s_i$. It is not an alcove walk in general, since some of the
$c_{i_\nu}$ will now carry the wrong exponents; the places where this will
happen are precisely those where the wall between $\mathbf a_{\nu-1}$ and
$\mathbf a_\nu$ has the form $H_{\alpha_i, k}$ for some $k\in\mathbb Z$ ---
we call those $\nu$ relevant.  Since $\langle \alpha_i, \lambda \rangle \ge
0$, all of those walls are crossed in the positive direction,
i.~e.~$\varepsilon_\nu = 1$ for all relevant $\nu$.
We want to write $c_i^- p_\lambda$ as a sum of walks. Let $\nu$ be the
minimal relevant index.  Since $c_{i_\nu}^+ = c_{i_\nu}^- + f_{i_\nu}^+$,
we have $c_i^- p_\lambda = p_1 + q_1$, where $p_1$ is obtained by replacing
$c_{i_\nu}^{+}$ in $c_i^- p_\lambda$ by
$c_{i_\nu}^{-}$, and $q_1$ is obtained by replacing it by
$f_{i_\nu}^{+}$. Then $q_1$ is indeed an alcove walk, and its
end point is $s_{\alpha_i+1}s_{\alpha_i}\epsilon^\lambda\mathbf a =
\epsilon^{\lambda-\alpha_i^\vee}\mathbf a$, and
since there is precisely one (positive) folding in $q_1$, in $\mathcal H$
the element $q_1$ is equal to $(q^{\frac{1}{2}} - q^{-\frac{1}{2}})
\theta_{\lambda - \alpha_i^\vee}$.  
On the other hand, $p_1$ will not be an alcove walk, in
general, and we have to repeat this procedure for the next relevant index.
Since $t_i = \Phi(c_i^-)^{-1}$, we have proved that 
\[ t_i^{-1} \theta_\lambda = c_i^- p_\lambda = (q^{\frac{1}{2}} -
q^{-\frac{1}{2}}) \theta_{\lambda - \alpha_i^\vee} + \cdots +
(q^{\frac{1}{2}} - q^{-\frac{1}{2}}) \theta_{\lambda - \langle
\alpha_i,\lambda\rangle\alpha_i^\vee} + \theta_{s_i\lambda}t_i^{-1}, \]
and since $t_i^{-1} = t_i + q^{-\frac{1}{2}}-q^{\frac{1}{2}}$, we get
\[ t_i \theta_\lambda = (q^{\frac{1}{2}} - q^{-\frac{1}{2}})\left(
\theta_{\lambda} + \theta_{\lambda - \alpha_i^\vee} + \cdots +
 \theta_{\lambda - (\langle
\alpha_i,\lambda\rangle-1)\alpha_i^\vee}\right) + \theta_{s_i\lambda}t_i, \]
which is the identity we had to prove.
\end{proof}

From the definitions it is clear that $t_{s_i} = \Phi(c_i^-)^{-1} =
\tilde{T}_i$ for $i=1, \dots, r$.  It is also clear that $\Theta_\lambda =
\theta_\lambda$ for $\lambda$ dominant or anti-dominant, and part (1) of
the proposition implies that $\Theta_\lambda = \theta_\lambda$ holds in
general.  Finally, it follows from part (4) that $\Phi(c_0^+) = \tilde{T}_0$.
Since the elements $\Theta_\lambda T_w$, $\lambda\in X_*$, $w\in W$ form a
basis of $\mathcal H$, we have obtained a proof of the Bernstein
presentation of $\mathcal H$.

\begin{remark} \label{gap}
The way of reasoning in \cite{R} is different from the above. Ram's
Proposition \cite{R} 3.2 (= Prop. \ref{Ram-Bernstein} above) 
shows (taking the Bernstein presentation of
the affine Hecke algebra as a piece of input), that the surjections from
the alcove walk algebra induce a (surjective) morphism $\mathcal H
\longrightarrow \widetilde{H}$.  It is implicitly stated that this is an
isomorphism, but the question whether this morphism is actually injective
is not addressed, as far as I can see.

Of course, the injectivity follows from Theorem \ref{indep} above.
Another way to approach this question is the following: 
It is clear that the quotient of $A$ by the
free $k$-submodule generated by the elements $f_i^+ -
(q^{\frac{1}{2}}-q^{-\frac{1}{2}})$, $i=0, \dots, r$ and $p-p'$ for $p$,
$p'$ non-folded alcove walks with ${\rm end}(p) = {\rm end}(p')$, admits a
basis of the form $\theta_\lambda t_w$, $\lambda\in X_*$, $w\in W$; hence
it would be enough to show that this $k$-submodule actually is an ideal.
(Again, it is clear from hindsight that this is true.)
\end{remark}

\subsection{The case of unequal parameters}
\label{unequalparams}

To keep the notation simple, we have considered only the case of equal
parameters. However, we can in a similar way as above consider a variant of
the alcove walk algebra with a parameter system. More precisely, assume
that a parameter system $L \colon \widetilde{W} \longrightarrow \mathbb Z_{\ge 0}$ is
given; we use the notation of Lusztig \cite{L}. Denote by $\mathcal H$ the
Hecke algebra associated with this system of parameters; see
loc.~cit. It is clear from the proof that Theorem \ref{indep} holds in this case,
too, when we use the elements $T_i$, and it is easy to see that as before
it holds with the $\tilde{T}_i=q^{-\frac{L(s_i)}{2}} T_i$, as well. 
The alcove walk algebra $A$ does not change, but we change 
the map $\Phi \colon A \longrightarrow \mathcal H$, namely we map
\[ c_i^+ \mapsto \tilde{T}_i,\quad  c_i^- \mapsto \tilde{T}_i^{-1},\quad
f_i^+ \mapsto q^{\frac{L(s_i)}{2}}-q^{-\frac{L(s_i)}{2}},\quad f_i^- \mapsto
q^{-\frac{L(s_i)}{2}}-q^{\frac{L(s_i)}{2}},\quad t_\tau\mapsto T_\tau.\] 
As above, $\Phi$ induces an
isomorphism $A/\mathscr J \cong \mathcal H$, where $\mathscr J \subseteq A$
is the two-sided ideal generated by the elements $f_i^+ -
(q^{\frac{L(s_i)}{2}}-q^{-\frac{L(s_i)}{2}})$, $i=0, \dots, r$, and $p-p'$
where $p$, $p'$ are non-folded alcove walks with the same end point.

It is slightly more difficult to prove the Bernstein relations in
$A/\mathscr J$ as in Proposition \ref{Ram-Bernstein} in this more general case in
terms of alcove walks. We have to prove that
\[ t_{s_i} \theta_\lambda = \left\{ 
\begin{array}{ll} \theta_{s_i\lambda}t_{s_i} +
(q^{\frac{L(s_i)}{2}}-q^{-\frac{L(s_i)}{2}}) \frac{\theta_\lambda -
\theta_{s_i\lambda}}{1-\theta_{-\alpha_i^\vee}} & \text{if } \alpha \not
\in 2X^*, \\
((q^{\frac{L(s_i)}{2}}-q^{-\frac{L(s_i)}{2}}) +
\theta_{-\alpha_i^\vee}(q^{\frac{L(\tilde{s}_i)}{2}}-q^{-\frac{L(\tilde{s}_i)}{2}}))\frac{\theta_\lambda -
\theta_{s_i\lambda}}{1-\theta_{-2\alpha_i^\vee}} & \text{if } \alpha 
\in 2X^* \end{array}
\right. \]
Roughly speaking, the same proof as above applies,
but obviously one has to be more careful in order to identify the
factors that come from the foldings: we need to know which $i_ \nu$ occur
for relevant $\nu$. The following two lemmas provide this information and
give a simple explanation for the appearance of the parameter associated to
$\tilde{s}$ when $\alpha \in 2 X^*$.

\begin{lemma}
Let $H \subseteq X_{*,\mathbb R}$ be an affine root hyperplane, 
and let $I \subseteq \{ 0, \dots, r \}$ be the set of types of faces 
with support $H$. Then $\{ s_i ;\ i \in I\}$ is a conjugacy class 
of simple reflections under the affine Weyl group $W_a$.
\end{lemma}

\begin{proof}
After applying an element of $W_a$ to $H$, if necessary, we may assume that
$H$ contains a face of the base alcove, say of type $i_1$. 

First assume that $F$ is a face with support $H$, of type $i_2$, say.
Choose $w\in W_a$ such that $F$ is a face of $w \mathbf a$. Since $F$ is
fixed by the reflection with respect to $H$, it is clear that the gallery
starting at $\mathbf a$ and consisting of crossing $H$, crossing the
faces of the types described by $w$, and finally crossing $H$ again, leads
to $w\mathbf a$. In other words, $s_{i_1} w s_{i_2} = w$, which shows that
$s_{i_1}$ and $s_{i_2}$ are conjugate.

On the other hand, assume that $s_{i_1} = w s_{i_2} w^{-1}$,
i.~e.~$s_{i_1} w = w s_{i_2}$. This shows that the image of $w\mathbf a$ 
under $s_{i_1}$ (which is just the reflection with respect to $H$) is
adjacent to $w\mathbf a$ by a face of type $i_2$. Since the alcoves are
adjacent, and $H$ lies between them, this face has to lie on $H$.
\end{proof}

The lemma shows that we can associate to each affine root hyperplane a
parameter $L(H) := L(s_i)$, where $i$ is the type of any face with support
$H$.

\begin{lemma}
Now suppose that $\alpha_i$ is a simple root.
\begin{enumerate}
\item If $\alpha_i \not\in 2X^*$, then for all $j\in \mathbb Z$, we have
$L(H_{\alpha_i,j}) = L(s_i)$.
\item If $\alpha_i \in 2 X^*$, then for $j$ even, $L(H_{\alpha,j})= L(s_i)$, 
and for $j$ odd, $L(H_{\alpha_i,j}) = L(\tilde{s}_i)$ (where we again use the
notation of \cite{L}).
\end{enumerate}
\end{lemma}

\begin{proof}
Consider the map $\mathbb Z \longrightarrow \mathbb Z_{\ge 0}$ mapping $j$ 
to $L(H_{\alpha_i,j})$. Since the types of faces are preserved by the action
of $W_a$, and since $s_{\alpha_i,j}(H_{\alpha_i,j-1}) = H_{\alpha_i,j+1}$,
the map factors through $\mathbb Z/2\mathbb Z$. Furthermore, the conjugacy
class corresponding to $H_{\alpha_i,0}$ is clearly the conjugacy class of
$s_i$, so $0$ maps to $L(s_i)$.

Now if $\alpha_i\not\in 2X^*$, then there exists $\lambda\in X_*$ with
$\langle \alpha_i,\lambda \rangle$ odd. Consider the translation of
$H_{\alpha_i,0}$ by $\lambda$. We have
\[ L(H_{\alpha_i,0}) = L(\epsilon^\lambda H_{\alpha_i,0}) =
L(H_{\alpha_i,1}). \]
If $\lambda$ is a coroot, such that translation by $\lambda$ is an element
of the affine Weyl group, then the previous lemma shows that the conjugacy
classes under $W_a$ associated with $H_{\alpha_i,0}$ and $\epsilon^\lambda 
H_{\alpha_i,0}$ coincide. In general, $\epsilon^\lambda \in \widetilde{W}$,
and the first equality follows from the fact that simple reflections
conjugate under $\widetilde{W}$ have the same parameter. The second
equality is clear since the pairing of $\alpha_i$ and $\lambda$ is odd.

Finally, we consider the case $\alpha_i\in 2 X^*$. This means that we are
in a very special case (cf.~\cite{L}, 2.4). The coxeter system $(W_a, S_a)$
is of type $\tilde{C}_r$, we may assume the root datum is adjoint, and $\tilde{s}_i$ is
the affine simple reflection $s_0$ (i.~e.~$\alpha_i=\alpha_r$ with the notation of
\cite{Bourbaki}).
Denoting by $\tilde{\alpha}$ the highest root, we have $\alpha_i +
\tilde{\alpha} = 2 \beta$ for a root $\beta$, and the image of
$H_{\alpha_i,-1}$ under the reflection $s_\beta$ is $H_{\tilde{\alpha},1}$.
Clearly, the conjugacy class associated with $H_{\tilde{\alpha}, 1}$ is the
conjugacy class of $s_0$, so the lemma is proved.
\end{proof}

\section{Minimal expressions}

\subsection{Existence of minimal expressions}

We say that $\Theta_\lambda \in \mathcal H$ has a \emph{minimal
expression}, if we can express it in the form
\[ \Theta_\lambda = \tilde{T}_{i_1}^{\varepsilon_1}\cdots
\tilde{T}_{i_\ell}^{\varepsilon_\ell} \tilde{T}_\tau, \]
where $\varepsilon_i \in \{ \pm 1 \}$ and where $s_{i_1}\cdots
s_{i_\ell}\tau$ is a reduced expression in $\widetilde{W}$; cf. the paper
\cite{HP} by Haines and Pettet, where the case of $\Theta_\lambda^-$
instead of $\Theta_\lambda$ is considered, and where minimal expressions
are interpreted in a sheaf-theoretic way, using Demazure resolutions of
Schubert varieties in the affine flag variety.  As an application of
Theorem \ref{indep} we get that minimal expressions always exist.

\begin{corollary} Let $\lambda$ be a coweight, and let $\epsilon^\lambda =
s_{i_1}\cdots s_{i_k}\tau$ be a reduced expression. Then 
\[ \Theta_\lambda = \tilde{T}_{i_1}^{\varepsilon}
\cdots \tilde{T}_{i_k}^\varepsilon \tilde{T}_\tau. \]
\end{corollary}

\begin{proof}
By definition
\[ \Theta_\lambda = \tilde{T}_{\lambda_1}
\tilde{T}_{\lambda_2}^{-1}, \]
where $\lambda = \lambda_1 - \lambda_2$, $\lambda_1$, $\lambda_2$ dominant,
and if we choose reduced expressions $t_{\lambda_1} = s_{j_1} \cdots
s_{j_\ell}\tau_1$, $t_{\lambda_2} = s_{j'_1} \cdots s_{j'_{\ell'}}\tau_2$, then it is
clear that 
\begin{eqnarray*} \tilde{T}_{\lambda_1}\tilde{T}_{\lambda_2}^{-1} & = &
\tilde{T}_{j_1} \cdots \tilde{T}_{j_\ell} \tilde{T}_{\tau_1}
\tilde{T}_{\tau_2}^{-1} \tilde{T}_{j'_{\ell'}}^{-1} \cdots \tilde{T}_{j'_1}^{-1} \\
& = & \tilde{T}_{j_1}^\varepsilon \cdots \tilde{T}_{j_\ell}^\varepsilon
\tilde{T}_{\tau(j'_{\ell'})}^{\varepsilon} \cdots
\tilde{T}_{\tau(j'_1)}^{\varepsilon} \tilde{T}_\tau.
\end{eqnarray*}
Now Theorem \ref{indep} immediately implies the result.
\end{proof}

It is clear that in an analogous way we can compute variants of the
$\Theta_\lambda$, where we write $\lambda$ as a difference not of dominant
coweights, but of coweights lying in some other fixed finite Weyl chamber.
In particular we see that the elements $\Theta_\lambda^-$ studied in
\cite{HP} (which are defined by $\Theta_\lambda^- =
\tilde{T}_{\lambda_1}\tilde{T}_{\lambda_2}^{-1}$, where $\lambda =
\lambda_1 - \lambda_2$ with $\lambda_1$, $\lambda_2$ anti-dominant) admit
minimal expressions as well.

\subsection{Applications}
The existence of minimal expressions has proved to be useful in several
occasions. For example, in \cite{GH1}, where the Jordan-H\"older series of
certain nearby cycles sheaves which arise naturally in the
Beilinson-Gaitsory deformation for the affine flag manifold to the affine
Grassmannian, were studied. The trace of Frobenius on the stalks of these
sheaves is a polynomial in $q$ and $q^{-1}$, and based on computational
evidence, Haines and the author conjectured that it is actually a
polynomial, and a sharp bound for its degree. In the
case where minimal expressions as above are available, we were able to
prove this conjecture; the results obtained above show that this proof is
valid in general. (Two quite different proofs of this fact are given in
\cite{GH2}.)

In the paper \cite{HP} by Haines and Pettet, there are several results
whose proofs rely on the existence of minimal expressions, as well.

\end{document}